\documentstyle[12pt]{article}

\newcounter{const}

\begin{document}

\begin{center}

{\bf Exponential Orlicz Spaces: New Norms and Applications.}\\
\vspace{3mm}

{\sc E.I.Ostrovsky.}\\
\vspace{3mm}
{\it Department of Mathematics and Computer Science, Ben Gurion University, \\
Israel,  Beer - Sheva, 84105, Ben Gurion street, 2, P.O. BOX 61,\\
E - mail: galaostr@cs.bgu.ac.il }\\
\vspace{3mm}

{\bf Abstract}
\end{center}
\vspace{3mm}

{\it The aim of this paper is investigating of Orlicz spaces with exponential  $ N \ $
 function and correspondence  Orlicz norm: we introduce some new equivalent  norms,
 obtain the  tail characterization,    
 study the product of functions in Orlicz spaces etc. \par
   We consider some applications: estimation of  operators in Orlicz 
 spaces and problem of martingales  convergence and divergence. } \\

\vspace{2mm}
\vspace{2mm}
 {\it Math. Sub. Classification} (2000): 47A45, 47A60, 47B10, 18D05. 60F10, 60G10, 60G50.\par

\vspace{2mm}
\begin{center}
{\sc 1. INTRODUCTION.} \par 
\end{center}
\vspace{2mm}
 Let $ (\Omega, F, {\bf P}) $ be a probability space. Introduce the following set 
of $ N \ - $  Orlicz functions:
$$
LW = \{N = N(u) = N(W,u) = \exp( W(\log u)) \},  \ u \ge e^2
$$
 where $ W $ is a continuous strictly increasing 
convex function in domain $ [ 2, \infty) $ such that 
$ u \to \infty \ \Rightarrow \lim_{u \to \infty} W_-^/(u) = \infty; $ here $ W_-^/(u) $
denotes the left derivative of the function $ W. $\par

 We define the function $ N(W(u)) $ arbitrary for the values $ u \in [0,e^2) $  
 but so that $ N(W(u)) $ will  be continuous convex strictly increasing and such that 
$$
u \to 0+ \ \Rightarrow N(W(u)) \sim C(W) u^2, 
$$
for some $ C(W) = const  \in (0,\infty). $ For $ u < 0 $ we define  as usually
$ N(W(u)) = N(W(|u|)). $ \par
 We denote  the set of all those $ N \ - $ functions as $ ENF: \ ENF = \{N(W(\cdot)) \} $
(Exponential N - Functions) and denote also the correspondence Orlicz space as   $ EOS(W) = $ 
Exponential Orlicz Space  or simple  $ EW $ with Orlicz (or, equally, Luxemburg) norm
$$
||\eta||L(N) = ||\eta||L(N(W)) = \inf_{v > 0} \{ v^{-1}(1+ {\bf E} N(W(\cdot), v \eta) \},
$$
where $ {\bf E}, {\bf D} $ denote the expectation and variance with respect to the 
probability measure $ {\bf P}: $ \par
$$
{\bf E} \eta = \int_{\Omega} \eta(\omega) \ {\bf P} (d \omega).
$$
 Let us introduce the following important function:
$$
\psi(p) = \psi(W,p)  = \exp(W^*(p)/p), \ \ p \ge 2, \eqno(1)
$$
where
$$
W^*(p) \stackrel{def}{=} \sup_{z \ge 2} (pz - W(z)) \ - 
$$
 is the Young  -  Fenchel transform of $ W. $  The function $ p \to p \ \log\psi(\cdot)= 
 W^*(\cdot) $ is continuous convex increasing and such that 
$$
p \uparrow \infty \ \Rightarrow p \log \psi(p) \uparrow \infty.
$$ 
 We denote the set of all those functions $ \{\psi\} = \{\psi(W,p) \} = \{ \exp(W^*(p)/p) \} $ 
by the symbol $ \Psi: $
$$ 
\Psi = \cup_{W} \{ \psi(W, p) \} = \{ \exp(W^*(p)/p) \}.
$$
 Inversely, $ W(\cdot) $ may be constructed by means of $ \psi: $
$$
W(p) = (p \log \psi(p))^* 
$$
(theorem of Fenchel - Moraux).\par
{\bf Definition}. {\it We introduce the so - called $ G(\psi) $ norm for arbitrary 
 function }  $ \psi(\cdot) \in \Psi  $ {\it and the correspondent $ G(\psi) $ space: 
$ \eta \in G(\psi) $ iff }

$$
||\eta||G(\psi) \stackrel{def}{=} \sup_{p \ge 2} |\eta|_p/\psi(p) < \infty,
$$
where $ |\cdot|_p $ denotes the classical $ L_p $ norm:
$$
|\eta|_p = {\bf E}^{1/p}|\eta|^p, \ p \ge 1.
$$
 In particular, $ |\eta|^2_2 = {\bf D} \eta + {\bf E}^2 \eta. $ It is easy to prove 
(see, for example, [3], p. 373; [12], p. 67) that $ G(\psi) $ is some (full) Banach space. \par
 Note that if there exists a family of measurable functions $ \{ \eta_{\alpha} \}, \ 
\alpha \in {\cal A} $ so that 
$$
\sup_{\alpha} |\eta_{\alpha}|_p < \infty,
$$
then there exists a $ G(\psi) \ - $ space, $ \psi \in \Psi $ such that $ \forall \alpha 
\ \Rightarrow \eta_{\alpha} \in G(\psi). $ For example, put
$$
\psi(p) = \sup_{\alpha} |\eta_{\alpha}|_p.
$$

{\it Remark 1.} In this paper the letters $ C, C_j $ will denote positive finite 
various constants which may differ from one formula to another and which do not 
depend on the essential parameters: $ x,u,p,z,\lambda. $\par

\begin{center}
\vspace{3mm}

{\sc 2. MAIN RESULTS.}\par
\end{center}
\vspace{3mm}

{\bf Theorem 1.} {\it The norms} $ ||\eta||L(N(W)) $ and $ ||\eta||G(\psi) $  {\it are equivalent. 
Further,} $ \eta \ne 0, \eta \in G(\psi) $  iff $ \exists C, C_1 \in (0,\infty) \ 
\Rightarrow \forall u > 2C $   
$$
{\bf P}(|\eta| > u) \le C_1 \exp(-W(\log (u/C))). \eqno(2)
$$
{\it Remark 2.} This result is a little  generalization of [3]; see also [16], 
p.305, as long as we do not suppose  that the function $ \eta = \eta(\omega), \ \omega 
\in \Omega $ to be exponential integrable:
$$
\exists \lambda > 0 \ \Rightarrow {\bf E} \exp(\lambda |\eta|) < \infty
$$
(so - called Kramer condition). \par
{\bf Proof.} a). Suppose  $ \eta \in EW, \eta \ne 0. $ Then for some $ C \in (0,\infty) \
\Rightarrow {\bf E} \exp(W (\log C |\eta|)) < \infty. $ Proposition (2) follows from Chebyshev 
inequality.\par
b). Inversely, let us assume that   $ \eta(\omega) $ is a measurable function such that
$$
{\bf P}(|\eta| > u) \le \exp(-W(\log u)), \ u \ge e^2.
$$ 
Then, by virtue of properties $ W(\cdot) $ we have:

$$
{\bf E} N(W(|\eta|/e^2)) = \int_{\Omega} N(W, |\eta(\omega)|/e^2) \ {\bf P} (d \omega) \le
$$
$$
 C_1 +  \sum_{k=2}^{\infty} \int_{e^k < |\eta| \le e^{k+1}} \exp((W(|\eta|/e^2)) \ 
 \ {\bf P}(d \omega) \le C_1+
$$

$$
 \sum_{k=2}^{\infty} \exp((W(k-1)) \ {\bf P}(|\eta|>k) \le C_2 + \sum_{k=2}^{\infty} 
\exp( W(k) - W(k+1)) < \infty.
$$
 Hence  $ ||\eta||L(N(W)) < \infty. $ \par 

c) Let $ \eta \in G(\psi); $ without loss of generality we can  assume  that 
$ ||\eta||G(\psi) = 1. $ We deduce: $ {\bf E}|\eta|^p \le \psi^p(p), $

$$
{\bf P}(|\eta| > x) \le x^{-p} \ \psi^p(p) = \exp (-p \log x + p \log \psi(p)), \ x > e^2,
$$
and, after the minimization over $ p: \ x \ge e^2 \ \Rightarrow $
$$
{\bf P}(|\eta|>x) \le \exp \left(-\sup_{p \ge 2}(p \log x - p \log \psi(p)) \right) = 
\exp(-W(\log x)).
$$
d) Let us next assume that  $ {\bf P}(|\eta|>x) \le \exp(-W(\log x)), x > e^2. $ We have:
$$
{\bf E}|\eta|^p \le C^p + \int_{e^2}^{\infty} p x^{p-1} \exp(-W(\log x) \ dx = 
$$

$$
 C^p + p \int_2^{\infty} \exp(py - W(y)) \ dy, \ p \ge 2.
$$
 Using  Laplace's  method and theorem of Fenchel - Moraux we get:

$$
{\bf E}|\eta|^p \le C^p + C_1^p \exp \left(\sup_{y \ge 2} (py - W(y)) \right)  = C^p+ 
$$

$$
 C_2^p \exp \left( W^*(p) \right) = C^p + C_2^p \exp(p \log \psi(p)) \le C_3^p \psi^p(p). \eqno(3)
$$
 Finally, $ ||\eta||G(\psi) \le C_3 < \infty.$ \par
{\it Remark 3.} If conversely 
$$
{\bf P}(|\eta| > x) \ge \exp(-W(\log x)), \ x \ge e^2,
$$
then for sufficiently large values of $p;\  p \ge p_0 = p_0(W) \ge 2 $ 
$$
|\eta|_p \ge C_0(W) \psi(p), \ C_0 \in (0,\infty).
$$
 
For arbitrary Orlicz spaces $ L(N), EOS(W) = EW, G(\psi) $ we denote correspondently 
$ L(N)^0, EW^0, G^0(\psi) $ a closure in $ L(N), EW, G(\psi)  $ norm 
the set of all bounded measurable functions.  It is known (see, for example, [15], p.75,
 [10], p.138) that in our conditions $ ({\bf P}(\Omega) = 1, \ N(W(u)) \sim 
C u^2 $ etc.) 
$$
LN(W)^0 = \{\eta: \forall k > 0 \ \Rightarrow {\bf E}N(W, \ |\eta|/k) < \infty. \}
$$ 

{\bf Theorem 2.} {\it Let} $ \psi \in \Psi. $ {\it We assert that} $ \eta \in EW^0, $ {\it 
or, equally,} $ \eta \in G^0(\psi) $ {\it if and only if} 

$$
\lim_{p \to \infty} |\eta|_p/\psi(p) = 0. \eqno(4)
$$
 {\bf Proof}. It is sufficient by virtue of theorem 1 to consider  only the case of the
$ G(\psi) $ spaces.\par
1. Denote $ GB^0(\psi) = \{ \eta: \lim_{p \to \infty} |\eta|_p/\psi(p) = 0.\} $
 Let  $ \eta \in  G^0(\psi), \eta \ne 0. $ Then for arbitrary 
$ \delta =  const > 0 $ there exists a constant $ K \in (0, \infty) $ such that 
$$
||\eta - \eta I(|\eta| \le K) \ ||G(\psi) \le \delta/2,
$$
 where for any event $ A \in F \ I(A) = 1, \omega \in A, \ I(A) = 0 $ if $ \omega \notin A. $ 
  Since $ | \eta I(|\eta| \le K )  | \le K, $ we deduce
$$
| \eta I(|\eta| \le K) |_p/\psi(p) \le K /\psi(p). 
$$
Using the  triangular inequality  we obtain for sufficiently large values 
$ p: \ p > p_0(\delta) = p_0(\delta,K): $

$$
|\eta|_p/\psi(p) \le \delta/2 + K/\psi(p) < \delta,
$$
as long as  $ \psi(p) \to \infty $ as $ p \to \infty.$  Therefore, $ G^0(\psi) 
\subset GB^0(\psi). $ \par
  (The set $ GB^0(\psi) $ is a closed subspace of $ G(\psi) $ with respect to the  
$ G(\psi) $ norm and contains all bounded random  variables). \par 
2. Conversely, assume  $ \eta \in GB^0(\psi). $ Let us denote  $ \eta(K) = \eta I(|\eta| > K), \
K \in (0,\infty). $ We deduce: 

$$ 
\forall Q \ge 2 \ \Rightarrow \ \lim_{K \to \infty} |\eta(K)|_Q = 0.
$$
Further,
$$
||\eta(K)||G(\psi) = \sup_{p \ge 2} |\eta(K)|_p/\psi(p) \le \max_{ p \in [2,Q]} 
|\eta(K)|_p/\psi(p) + 
$$

$$
 \sup_{ p > Q} |\eta(K)|_p/\psi(p) \stackrel{def}{=} \sigma_1 + \sigma_2;
$$
 
$$
\sigma_2 = \sup_{p > Q} |\eta(K)|_p /\psi(p) \le \sup_{p \ge Q} \left( |\eta|_p/\psi(p) \right)
 \le \delta/2 
$$
for sufficiently large $ Q $ as long as $ \eta \in GB^0(\psi). $  Further,
$$
\sigma_1 \le \max_{p \in [2, Q]} |\eta(K)|_p/\psi(2) \le |\eta(K)|_Q/\psi(2) \le 
\delta/2
$$
for sufficiently large $ K = K(Q). $ Therefore, 

$$
\lim_{K \to \infty} ||\eta(K)||G(\psi) = 0, \ \ \Rightarrow \eta \in G^0(\psi).
$$
 Hence $ GB^0(\psi) \subset G^0(\psi). $ \par 
 Let now $ \{\eta_a \}, \ a \in {\cal A} $ be some {\it family } of functions from the
$ G^0(\psi) $ space.\par
{\bf Theorem 3.} {\it Let } $ \psi \in \Psi.$
{\it In order to a family} $ \{ \eta_a \} $ {\it of a function belonging 
to the} $ LG(\psi) $  {\it space  has the  uniform absolute 
continuous norm in this  space,  briefly: }
$ \{ \eta_a \} \in UCN(G(\psi)), $  {\it it is necessary and sufficient:}

$$
\lim_{p \to \infty}  \sup_{a \in {\cal A} } |\eta_a|_p /\psi(p) = 0.\eqno(5)
$$
 
{\bf Proof.} Recall that by definition  $ \{\eta_a \} \in UCN(G(\psi)) $  if

$$
\lim_{\delta \to 0+} \sup_{V: {\bf P}(V) < \delta} \sup_a ||\eta_a I(V)||G(\psi) = 0.
$$

1. Let the condition (5) be satisfied, then there exists a function 
$ \epsilon = \epsilon(p) \to 0 $ as $ p \to \infty $  such that $ \forall a \in {\cal A} $  and 
$ \forall p \ge 2 \ \Rightarrow |\eta_a|_p \le \epsilon(p) \psi(p). $ It follows that for all 
$ Q \ge 2 $  the family of functions  $ |\eta_a|^Q $ is uniform integrable.  Let $ V $ be 
an arbitrary measurable set: $ V \in F $ with sufficiently small measure: $ {\bf P}
(V) \le \delta, \ \delta \in (0,1/2). $ We have:

$$
\sup_a ||\eta_a I(V)||G(\psi) \le \sup_a \max_{p \le Q} |\eta_a|_p/\psi(p) + 
$$

$$
\sup_a \sup_{p > Q}|\eta_a|_p/\psi(p) \stackrel{def}{=} s_1 + s_2;
$$

$$
s_2 \le \sup_{ p \ge Q} \epsilon(p) \to 0, \ Q \to \infty;
$$

$$
s_1 \le \sup_a \ |\eta I(V)|_Q /\psi(2) \to 0, \ \delta \to 0+.
$$
 Hence the family $ \{\eta_a \} $ has the  uniform absolute continuous 
norms in our Orlicz space $ G(\psi). $\par  
 2. Assume now  the family $ \{ \eta_a \}, \ \eta_a \in G^0(\psi) $
belongs to $ UCN(G(\psi)): \ \{\eta_a \} \in UCN(G(\psi)). $
 Then this family is uniformly finitely approximate in $ G(\psi) $  norm: 
$$
\lim_{K \to \infty} \sup_a ||\eta_a \ I(|\eta_a| > K) \ ||G(\psi) = 0.
$$
 Further proof is the same  as in the theorem 3. \par
 Recall here the definition for  two $ N - $ functions $ M = M(u), \  N = N(u) $
(see, for example, [8], chapter 2, section 13, [9], p.144): the function $ N(\cdot) $ is 
called essentially greater than $ M(\cdot): N = N(\cdot) >> M = M(\cdot) $ or equally
$ M(\cdot) $ decreases essentially more rapidly then $ N(\cdot): \  M << N, $ if
$$
\forall \lambda > 0 \ \Rightarrow \lim_{u \to \infty} M(\lambda u)/N(u) = 0.
$$

 {\bf Theorem 4.} {\it Let $ \psi(\cdot) = \psi_N(\cdot), \ \nu(\cdot) = \nu_M(\cdot) $ be 
 two functions of the 
classes $ \Psi $ with correspondent $ N \ - $ Orlicz functions $ N(\cdot), M(\cdot): $

$$
N(u)= N_{\psi}(u) = \exp \{[p \log \psi(p)]^*(\log u)\}, \ 
$$
$$
M(u) = N_{\nu}(u)  = \exp \{[p \log \nu(p)]^*(\log u) \}.
$$
 We assert that}  $ \lim_{p \to \infty} \psi(p)/\nu(p) = 0 $ {\it if and only if} $ N >> M. $ \par
{\bf Proof.} 1. Assume that  $ \lim_{p \to \infty}\psi(p)/\nu(p) = 0. $ 
 Denote $ \epsilon(p) = \psi(p)/\nu(p), $ then $ \epsilon(p) \to 0, \
p \to \infty. $ \par
 Let $ \{\eta_a \}, a \in {\cal A} $ be arbitrary  bounded in the $ G(\psi) $ sense 
set of functions:
$$
\sup_a ||\eta_a||G(\psi) = \sup_a \sup_{p \ge 2}|\eta_a|/\psi(p) =  C < \infty, 
$$
then 

$$
\sup_a |\eta_a|_p/\nu(p) \le C \ \epsilon(p) \to 0, \ p \to \infty.
$$
 It follows from theorem 3 that $ \forall a \in {\cal A} \ \ \eta_a \in G^0(\nu) $ and that 
the family $ \{\eta_a \} $ in the space $ G^0(\nu) $  has the uniform absolute continuos 
norm: $ \{\eta_a \} \in UCN(G(\nu)).$ Our assertion it follows from lemma 13.3 in the book 
[8].\par
2). Inversely, let $ N >> M. $ Let us consider the measurable function $ \eta: 
\Omega \to R^1 $ such that $ \forall x \ge C $ 

$$
C_1 \exp(- C_2[p \log \psi(p)]^*(\log x)) \le {\bf P}(|\eta|> x) \le 
$$

$$
 C_3 \exp(-C_4[p \log \psi(p)]^*(\log x)).
$$
 It follows from theorem 1 that $ \eta \in G(\psi) $ and $ C_5 \psi(p) \le |\eta|_p \le 
C_6 \psi(p), \ p \ge 2. $ Since $ ||\eta||G(\psi) < \infty, \  M << N, $ we deduce that 
$ \eta \in G^0(\nu) $ ([8], theorem 13.4). It follows from theorem 2 that 
$$  
 \lim_{p \to \infty} |\eta|_p/\nu(p) = 0.
$$
Thus $ \ \ \lim_{p \to \infty} \psi(p) /\nu(p) = 0. $ \par
 Obviously, if for some $ C_1, C_2: 0 < C_1 \le C_2 < \infty $  and for all $ p \ge 2 $

$$
C_1 \le \psi(p)/\nu(p) \le C_2,
$$
 then the norms $ ||\cdot||G(\psi) $ and $ ||\cdot||G(\nu), $ or, equally, $ ||\cdot||L(N_{\psi}) $
and $ ||\cdot||L(N_{\nu}) $  are equivalent. \par

{\bf Theorem 5.} {\it Let $ N \in LW $ and $ \eta $ be a random variable (measurable function
$ \eta: \Omega \to R^1) $ such that for sufficiently large values } $ u: \ u \ge C_0 \ \Rightarrow $
$$  C_1 \exp(-W(\log u/C)) 
\le {\bf P}(|\eta| > u) \le C_2 \exp(-W(\log u/C_3)), \eqno(6)
$$
{\it then } $ \eta \in EW \setminus EW^0. $ \par
{\bf Proof.}  It follows from the right - side of inequality (6) in accordance to the theorem 1 
that $ \eta \in EW = G(\psi). $ Further, as well as in the proof of theorem 4 we obtain

$$
\overline{\lim}_{p \to \infty} |\eta|_p/\psi(p) > 0.
$$
 It follows from theorem 2 that $ \eta \notin G^0(\psi). $ \par
 Let us consider some concrete examples. The so - called $ EL_m, m = const > 0 $ spaces are
very important subclasses of the $ EW $ spaces.
 By definition of the $ EL_m $ spaces,
$$
W(\log u) = W_{L,m}(\log u) \stackrel{def}{=} u^m L(u), \ u \ge e^2,
$$
and correspondently $ N(u) = \exp (u^m L(u)), $
 where $ L(\cdot) $ is a slowly varying at $ u \to \infty $ continuous function such that 
$ L(0) > 0, \ u^m L(u) \uparrow \infty $ at $ u \uparrow \infty $ and 
$$
\lim_{u \to \infty} L(u/L(u))/L(u) = 1. \eqno(7)
$$
 We denote the set  of all those functions as $ SV_m, \ SV \stackrel{def}{=}
\cup_{L} \{L(\cdot) \}. $ \par
  It is known ([12], p. 25) that in the case $ m > 1, \ L \in SV_m $ 
the centered random variable  $ \eta, \ {\bf E}\eta = 0 $ belongs to the space $ L_m:
\eta \in EL_m $ if and only if  $ \exists C \in (0,\infty),  \forall \ \lambda \ge 2 $

$$
 {\bf E} \exp( \pm  \lambda \eta ) \le \exp \left( C \lambda^{m/(m-1)} \ 
L^{-1/(m-1)}(\lambda^{1/(m-1)}) \right). \eqno(8)
$$

 The following particular cases of $ EL_m $ spaces are very convenient in the practical using.
 Define for $ p,u \ge 2, m > 0, r \in R^1 $

$$
\psi_{m,r} = \psi_{m,r}(p) = p^{1/m} \log^r p, \ \  G_{m,r}= G(\psi_{m,r}), G_m = G_{m,0},
$$

$$
||\eta||_{m,r} = ||\eta||G(\psi_{m,r}), \ \ ||\eta||_m = ||\eta||_{m,0},
$$

$$
N_{m,r}(u) = \exp \left( u^m \ (\log(C(m,r) + u))^{-mr} \right).
$$
 The Orlicz spaces $ (\Omega,F,{\bf P}; L(N_{m,r}(\cdot)) $ and $ G_{m,r} $ are isomorphic and 
the correspondence norms are equivalent

$$
C_1(m,r) ||\eta||_{m,r} \le ||\eta||L(N_{m,r}) \le C_2(m,r)||\eta||_{m,r}.
$$ 

 Moreover, $ \eta \in G_{m,r} $ iff  $ \forall x \ge 2 \ $ 
$$ 
 {\bf P}(|\eta| > x) \le \exp \left( -C_3(m,r) x^m (\log(C(m,r) + x))^{-mr} \right).
$$

 Another example of  $ EOS(W) $ spaces are so - called $ V(Z, \beta) $ spaces: 
$ V(Z,\beta) = G \left(\psi^{(Z,\beta)} \right), $ where by definition

$$
\psi^{(Z,\beta)}(p) = \exp \left(Z p^{\beta} \right), \ Z,\beta = const > 0.
$$
 From theorem 1 follows that $ \eta \in V(Z,\beta), \eta \ne 0 $ iff $ \exists C \in 
(0,\infty) \ \Rightarrow \ \forall x \ge 2C $

$$
{\bf P}(|\eta|>x) \le 2 \exp \left( - Z^{-1/\beta} (1+\beta)^{1+1/\beta} 
(\log x/C)^{1+1/\beta} \right).
$$ 

 We continue investigating  the properties of $ EL_m $ spaces. Let $ \epsilon(k) $ be a
Rademacher sequence, i.e. the sequence of independent random variables with distributions:

$$
{\bf P}(\epsilon(k) = 1) = {\bf P}(\epsilon(k) = -1) = 1/2.
$$ 
 
 Let also $ B = const \in (0.5; 1) $ and  $ L = L(u) $ be a function belonging  to
the class $ SV: L(\cdot) \in SV. $ Introduce the random variable $  \xi $ by the formula

$$
\xi = \sum_{k=2}^{\infty} k^{-B} \ L(k) \ \epsilon(k). \eqno(9)
$$

 Denote
$$
  {\tilde L} (u) = L^{-1/(1-B)} \left(u^{1/(1-B)} \right). 
$$

{\bf Theorem 6.}  {\it There exist $ C_1, C_2 \in (0,\infty), C_1 \le C_1 $ so that } 
$ \forall u \ge 2 $

$$
\exp \left(-C_2 u^{1/(1-B)} {\tilde L}(u) \right) \le {\bf P}(|\xi| \ge u) \le 
$$

$$
 \exp \left( - C_1 u^{1/(1-B)} {\tilde L}(u) \right). \eqno(10)
$$
{\bf Proof.} Since

$$
\sum_{k=2}^{\infty} {\bf D} \left(k^{-B} \ L(k) \ \epsilon(k) \right) < \infty,
$$
then  there exists the r.v. $ \xi $ and has a symmetric distribution. \par
 Introduce for all values $ \lambda \in R^1 $ the function

$$
\varphi(\lambda) \stackrel{def}{=} \log {\bf E} \exp(\lambda \xi) = \sum
_{k=2}^{\infty} \log \cosh \left( k^{-B} L(k) \lambda \right). 
$$
 We have at $ \lambda \to + \infty: \ \varphi(\lambda) \sim $

$$
  \int_2^{\infty} \log \cosh \left( x^{-B} L(x) \lambda  \right) \ dx \sim 
 \lambda^{1/B}  \int_0^{\infty} \log \cosh \left(z^{-B} L \left(z \lambda^{1/B}\right)
\right) \ dz \sim 
$$

$$
 \lambda^{1/B} \int_0^{\infty} \log \cosh \left(z^{-B} L \left(\lambda^{1/B}) \right) 
\right) dz = \lambda^{1/B} L^{1/B} \left(\lambda^{1/B} \right) \ C(B), 
$$
where 
$$
C(B) = \int_0^{\infty} \log \cosh \left( z^{-B}\right) \ dz \in (0,\infty).
$$
since $ B \in (0.5; \ 1). $ \par
  Assertion (10) follows from the main result of paper [2]. \par

{\bf Theorem 7.} {\it Let $ \xi, \eta $ be a random variables belonging to the $ EL_m $
space, $ m > 0, \ L(\cdot) \in SL.$ Denote $ \tau = \xi \eta. $ \\
(A). We assert that } for $ x \ge 2 $

$$
 \ {\bf P}(|\tau| > x) \le \exp \left( - C x^{m/2} L(\sqrt{x}) \right). \eqno(11)
$$
(B). {\it Inversely, assume that the random variables $ \xi, \ \eta $ are independent,
identically symmetrical distributed and $ \exists \ L(\cdot) \in SV, \ \forall x \ge 2 
 \ \Rightarrow $
$$
\exp \left(- C_2 \ x^m L(x)  \right) \le {\bf P}(|\xi| > x) =
$$

$$
 {\bf P}(|\eta| > x) \le \exp \left(- C_1 \ x^m L(x) \right), 
$$
where $ 0 < C_1 < C_2 < \infty $ (the case $ C_1 = C_2 $ is trivial). Statement:}
  $ x \ge 2 \ \Rightarrow $

$$
 \ {\bf P}(|\tau| > x) \ge \exp \left( - C_3 x^{m/2} \ L(\sqrt{x}) \right). \eqno(12)
$$
{\bf Proof.} As long as 
$$
{\bf P}(\xi^2 > x) \le  \exp \left(-C_1 x^{m/2} L(\sqrt{x}) \right),
$$ 
 the random variables  $ \xi^2, \ \eta^2 $ belong to the Orlicz space $ EL_{m/2} $ with $ L \ $ 
function $ L(\sqrt{x}), \ x \ge 2. $  The first assertion (A) it follows from the linear 
properties of the $ EL_{m/2} $ spaces and from the elementary relation:
  $ 2 \xi \eta = (\xi + \eta)^2 - \xi^2 - \eta^2. $ \par  
  Now we are going to prove the assertion (B). Let the value $ x $ be sufficiently large. We 
have for some $ Y = const  \in \left(1, \sqrt[m]{C_2/C_1} \right) $ and $ k = 2,3,\ldots, $ using 
the full probability formula, since the function $ L(\cdot) $ is slowly varying: 
$ \overline{P} := {\bf P}(|\tau|> x) \ge $
$$
 \sum_k {\bf P}((|\xi \eta| > x) /(Y^k \le |\eta| < Y^{k+1}) 
\times {\bf P}(|\eta| \in \left[Y^k, Y^{k+1} \right) \ge  
$$

$$
 \sum_k {\bf P} \left(|\xi| > x Y^{-k-1} \right) \times {\bf P} 
\left( |\eta|\in \left[Y^k, Y^{k+1} \right) \right) \ge
$$

$$
 \sum_k \exp \left(-C_2 Y^{-m} (x/Y^k)^m L(x/Y^{k+1}) \right) \times
$$

$$
\times \left[ \exp \left(-C_2 Y^k \right) - \exp \left(-C_1 Y^{k+1} L(Y^{k+1} 
\right) \right] \ge
$$

$$
 \sum_k \exp \left( -C_4 (x/Y^k)^m L(x/Y^k) - C_5 V^k L(Y^k) \right)
$$
 Choosing in this sum the member with $ k =  k_1(x) = Ent[\log_Y(\sqrt{x})], $ 
where $ Ent[z] $ is an integer part of $ z, $  
 we obtain for sufficiently large values $ x $ 

$$
\overline{P} \ge \exp \left(- C_6 x^{m/2} L(\sqrt{x}) \right).  
$$

\vspace{2mm}
\begin{center}
{\sc 3. APPLICATIONS TO THE OPERATOR'S THEORY.} \par
\end{center}
\vspace{2mm}
 Let $ (\Omega_1, F_1, \mu_1) $ and $ (\Omega_1, F_2, \mu_2) $ be two probability 
spaces, $ L^j_p = L_p(\Omega_j), \ G^j_m = G_m(\Omega_j), j = 1,2; $
$$
|f|^{(j)}_p = |f|L^j_p, \ \  ||f||^{(j)}_m = ||f||G^{j}_m,
$$
 and $ Q = Q[f], \ f: \Omega_1 \to R^1, \ Q[f]: \Omega_2 \to R^1  $ be an operator, not 
necessary linear, defined on the set 
$$
Dom (Q) = \cap_{p \ge 2}L^1_p 
$$
 with image 
$$
Im(Q) = \{Q[f], \ f \in Dom(Q) \}  \subset  \cup_{p \ge 2} L^2_p = L_2^2.
$$

{\bf Theorem 8.}  {\it Suppose that there exist some  constants $ a \ge 0, \  b,d, C, C_1 
\in (0,\infty) $ such that } $ \forall p \ge p_0 = const \ge 2 \ \Rightarrow $

$$
|\ Q[f] \ |^{(2)}_p \le C_1 \ p^a \ \left[| \ f \ |^{(1)}_{C p^b} \right]^d. \eqno(13)
$$
 {\it Then for all values $ m > 0 $ the operator $ Q $ may be defined on the Orlicz 
space $ G_m^1 $  with  values into Orlicz space $ G_n^2, \ n = m/(am + bd) $ and }
$ \forall f \in G_m^1 \ \Rightarrow $

$$
||Qf||^{(2)}_n \le C \left[ ||f||^{(1)}_m \right]^d. \eqno(14)
$$
 {\bf Proof} is very simple. Let $ f \in G^1_m, \ ||f||^{(1)}_m = 1, $ then 
$ \forall p \ge 2 \ \Rightarrow |f|^1_p \le C_2 p^{1/m}. $  It follows from condition 
(13) that 
$$
| \ Qf |_p^{(2)} \le C_3 p^a \cdot \left( p^{bd/m} \right) \left[||f||^{(1)} \right]^d_m 
= C_3 p^{1/n}.
$$
Our statement (14)  follows from theorem 1. \par
 There are many operators  satisfying  the 
 condition (13), for example: Hilbert operator ([5], p.119),  singular 
integral operators of a type Kalderon - Zygmund, Hardy operator  ([18], p.  42), operator of 
solution of linear or non - linear evolution equation [6] and  many others. In 
all known cases the values $ a $ are only $ a= 0; a = 1/2 $ and $ a = 1. $ \par 
 A very simple example of operator $ Q[\cdot] $ with 
 $ d \ne 1: \ Q[f](\omega) = f^d(\omega).$ \par
 Let us consider now in detail partial Fourier sums in the bounded domain
$ [0; 1] $ with  usually Lebesque measure $ \mu: $

$$
S_N[f](x) = 0.5 a_0 + \sum_{k=1}^N (a_k \cos (2 \pi x) + b_k \sin( 2 \pi x)),
$$
where $ a_k, b_k $ are Fourier's coefficients of the (integrable) function $ f. $
It is  well - known (theorem of M. Riesz): 

$$
p \ge 2 \ \Rightarrow | \ H[f] \ |_p \le C \ p \ |f|_p,
$$
where the symbol $ H[\cdot] $ denotes the Hilbert transform for a functions  defined on 
the interval $ (0,1) $ and, equally, (see, for example,  [5], p. 119 - 121)  for some 
other $ C \in (0,\infty) $ 

$$
\sup_{N \ge 1} |S_N[f]|_p \le C \ p \ |f|_p,
$$ 
where $ C $ is an {\it absolute} constant, i.e. here $ a = 1 $ uniformly on $ N. $
 It is proved in the paper [14]
that this estimate is exact at $ p \to \infty. $\par
 
By virtue of theorem 8 we conclude: $ \forall m > 0 $ and $ \forall f(\cdot) 
\in G_m \  \Rightarrow $

$$
|| \ H[f] \ ||_{m/(m+1)} \le C_0(m) \ ||f||_m, 
$$

$$
\sup_{N \ge 1} || \ S_N[f] \ ||_{m/(m+1)} \le C_1(m) ||f||_m,
$$
{\bf Lemma 1.}  {\it The previos constant $ m/(m+1) $ is optimal in the case
$ m \ge 1.$  In detail, } $ \forall \  m \ge 1 \ \exists g \in G_m \ \Rightarrow 
\forall \Delta > 0 $

$$
|| H[g]||_{(m+\Delta)/(m+1)} = \infty.
$$

{\bf Proof.} Let us introduce the function 
$$
g(x) = g_m(x) = \sqrt[m]{|\log x |}.
$$
 Since 
$$
\forall u \ge 0 \ \Rightarrow \mu \{x: g_m(x) > u \} = \exp \left(- C u^m \right),
$$
we conclude: $ g_m(\cdot) \in G_m \setminus  G_m^0 $ (theorem 5). Further,
 it is very simple to verify using the explicit view of Hilbert transform  ([5], p. 
112) that 
$$
C_1 \left(|\log x |^{(m+1)/m} + 1 \right) \le   |H[g_m](x)| \le 
$$
$$
C_2 \left( |\log  x|^{(m+1)/m} + 1 \right).
$$
 Hence $ \forall u \ge 2 $ 

$$
\exp \left( - C_3 u^{m/(m+1)} \right) \le  \mu \{ x, H[g_m](x) > u \} \le \exp 
\left( -C_4 u^{m/(m+1)}  \right)
$$
and again follows from theorem 5: 
$$
H[g_m] \in G_{m/(m+1)} \setminus  G^0_{m/(m+1)}.
$$
Hence $ H[g_m](\cdot) $ does not belongs to the space $ G_{(m+\Delta)/(m+1)} $ for all 
values $ \Delta \in (0,\infty).$ \par
 Consequently, the Fourier sums $ S_N[f] $ converge to $ f $ for all functions $ f \in G_m $
in the {\it other } Orlicz  space $ G(\nu): $ 
$$
\lim_{N \to \infty} ||S_N[f] - f||G(\nu) = 0 
$$
if for example 
$$ 
\lim_{p \to \infty} \nu(p) \ p^{-1/(1+1/m)} = \infty 
$$
 Now we construct the examples  (for all values $ m > 0) $ of the  function  $  g \in G_m $ 
such that  

$$ 
\lim_{N \to \infty} ||S_N[g] - g||_m \ne 0. \eqno(15) 
$$
Put again $ g(x) = g_m(x), $ but here $ m \in (0, \infty). $ It is evident that 
$ g(\cdot) \in G_m \setminus G^0_m $ (theorem 5). As long as the trigonometrical system is 
bounded, the assertion (15) is true. \par

\vspace{2mm}
\begin{center}
{\sc 4. APPLICATIONS TO THE THEORY OF MARTINGALES.} \par
\end{center}
\vspace{3mm} 

 Let $ (\Omega,F, {\bf P}) $ be an {\it arbitrary}  probability space with some filtration 
(flow of $ \sigma - $ algebras) $ (\emptyset, \Omega) =F_0 \subset F_1 \subset F_2 \ldots 
\subset F_n  \ldots \subset F $
and $ (S_n. F_n) $ be a centered  martingale: $ \forall n =1,2,\ldots \ {\bf E} S_n = 0, 
 \ {\bf E} |S_n| < \infty, $ and 

$$
{\bf E} S_n/F_{n-1} = S_{n-1}\ \ (mod \ {\bf P}), \ n=1,2,\ldots.
$$

  It is well known ([7], p.18; [13] ) that if for some $ p > 1 \ 
\sup_n |S_n|_p < \infty $ then the  martingale $ S_n $ converges a.e. and in $ L_p $ norm:
$ \exists S = \lim  S_n \ (mod \ {\bf P }), \ $ at $ n \to \infty; \ S \in L_p $ and

$$
\lim_{n \to \infty}  |S_n - S|_p \to 0.
$$
 
Some generalizations of this statement are obtained in  the publications [13], [11], p. 217;
in particular, on the Orlicz spaces $ (\Omega,F, {\bf P}, N(\cdot)) $
with $ N \ $ function (convex, even etc.)  belonging to the  $ \Delta_2 \cap \nabla_2 $ class:
$$
\Delta_2 = \{N: \exists \ (u_0 > 0, \beta < \infty), \  \forall u \ge u_0 \Rightarrow N(2u) \le 
\beta N(u) \},
$$
$$
\nabla_2 =  \{N: \exists \ (u_0>0, l > 1), \ \forall u \ge u_0 \ \Rightarrow N(u) \le N(lu)/(2l) \}.
$$
For example, if $ N(u) = N_m(L), \ m > 1, $ then $ N(\cdot) \in \Delta_2 \cap \nabla_2. $\par
Recall that if $ N(\cdot) \in \Delta_2  \cap \nabla_2 $ and the measure $ {\bf P } $ is diffuse 
that the Orlicz space $ L(N) $ is separable and reflexive. \par
 (Note that the considered above Exponential Orlicz $ N \ - $ functions $ N \in LW $
 do not satisfy the $ \Delta_2 $ condition.) \par 
 It is proved, more exactly,  in [13], [11], p. 217 that if $ N \in \Delta_2 \cap \nabla_2, $  
then for  arbitrary  martingale $ (S_n,F_n) $
$$
\sup_n ||S_n||L(N) < \infty  \ \Rightarrow \lim_{n \to \infty} ||S_n - S||L(N) = 0. 
$$

 We investigate in this section  the convergence $ S_n $ to $ S $ in the Orlicz  
spaces with $ N - $ function $ M(\cdot) $ {\it without $ \Delta_2 $ condition.} Namely,
we describe for any fixed  exponential $ N - $ function $ N = N(\cdot) \in ENF $ the 
set of  {\it other}  $ N - $ functions $ M = M(u) $ such that for {\it arbitrary} 
martingale $ (S_n,F_n) $ holds the following implication: 
$$
\sup_n ||S_n||L(N) < \infty  \ \Rightarrow \lim_{n \to \infty} ||S_n - S||L(M) = 0. \eqno(16)
$$
 Generally speaking, in the considered case $ N \in ENF \ \Rightarrow 
\lim_{u \to \infty} M(u)/N(u) = 0. $
Further  we construct some examples when  functions $ M,N $ are not equivalent.\par 
 
 It follows  from the theory of Orlicz spaces ([8], theorem 13.4)  that the condition 
$ M(\cdot) << N(\cdot) $ is sufficient for the implication (16).
 Further we find some other conditions (necessary conditions and sufficient conditions).\par

 For some $ \psi \ - $ function $ \psi \in \Psi $ mentioned above, $ p \ge 2, \ 
\delta \in (0,1) $ we define a new function $ R = R(\delta,p,\psi(\cdot)) = $ 
$$
\inf \left[ \delta^{2/(p \beta + 2)} \ \psi^{p \beta/(p \beta + 2)}
(\alpha p): \ \alpha, \beta > 1, 1/\alpha + 1/\beta = 1 \right],
$$
where {\it inf} is calculated over all values $ (\alpha, \beta)$ such that 
$ \ \alpha,\beta > 1, \ 1/\alpha + 1/\beta = 1. $ \par
{\bf Theorem 9. } {\it Let $ \nu = \nu(p) $ be some function, $ \nu \in \Psi, $ 
so that}
$$
\lim_{\delta \to 0+} \sup_{p \ge 2} R(\delta,p,\psi(\cdot))/\nu(p) = 0. \eqno(17)
$$
{\it Then for all martingales $ (S_n, F_n) $ such that }
$$
\sup_n ||S_n||G(\psi) < \infty \eqno(18)
$$
follows the implication:
$$ 
n \to \infty \ \Rightarrow \ || S_n - S|| G(\nu) \to 0, \eqno(19)
$$
 {\it i.e. the martingale $ S_n $ converges in the sense of $ G(\nu) $ norm.} \par
{\bf Proof.} From our condition (and classical theorem of Doob) follows that 
there exists  $  S = \lim_{n \to \infty}S_n $ a.e. and $ \sigma_n^2 \to \sigma^2, $ where 
$ \sigma^2_n = {\bf D}S_n \le  \sigma^2 = {\bf D} S.$ Let us denote $ \gamma^2 = \gamma^2_n =
\sigma^2 - \sigma^2_n = |S-S_n|^2_2, \ (\gamma_n  \to 0, \ n \to \infty); $

$$
K = \sup_n ||S_n||G(\psi) < \infty.
$$
 
 We obtain using Chebyshev inequality:  $ \ \varepsilon > 0 \ \Rightarrow $
$$
{\bf P}(|S_n-S|>\varepsilon) \le (\sigma^2-\sigma^2_n)/\varepsilon^2 = 
\gamma_n^2/\varepsilon^2.
$$ 
 Put $ \xi_n = S_n - S. $ We have by virtue of Doob's inequality, 
since $ \sup_n ||S_n||G(\psi) < \infty: $
$$
|| \max_{l \le n} \  S_l \ ||G(\psi) \le \sup_{p \ge 2} \{ [p/(p-1)] \ 
|S_n|_p/\psi(p) \}  \le 
$$

$$
 2 \sup_{p \ge 2} |S_n|_p/\psi(p) = 2 || \ S_n \ ||G(\psi) \le 2K,
$$
hence
$$
|| S||G(\psi) \le  || \sup_n |S_n| \ ||G(\psi) \le 4K,
$$
$$
\sup_n ||\xi_n||G(\psi) \le \sup_n || S_n||G(\psi) + ||S||G(\psi)  \le 5K.
$$
 Further we have for all $ p \ge 2, \varepsilon > 0 $ 
$$
{\bf E} |\xi_n|^p = \int_{\Omega} |\xi_n|^p d {\bf P} \le \varepsilon^p + 
\int _{\Omega} |\xi_n|^p \ I( |\xi_n|\ge \varepsilon) \ d {\bf P}.
$$
 We get estimating the right side by the $ H\ddot{o}lder $ {\it inequality} 
$$
\int_{\Omega} |\xi_n|^p I(|\xi_n| > \varepsilon) \ d {\bf P} \le 
\left( \int_{\Omega} |\xi_n|^{\alpha p} \ d {\bf P}  \right)^{1/\alpha} \times
$$
$$
\times \left({\bf P}(|\xi_n| > \varepsilon   \right)^{1/\beta} \le |\xi_n|^p_{\alpha \ p} \
(\gamma_n/\varepsilon)^{2/\beta} \le 
 5^p \ K^p  \ \psi^p(\alpha p) \ (\gamma_n/\varepsilon)^{2/\beta},
$$
where as above $ \alpha,\beta > 1, 1/\alpha + 1/\beta = 1. $ Therefore 
$$
{\bf E} |\xi_n|^p \le \varepsilon^p + 5^p \ K^p \ \psi^p(\alpha p) \ \gamma_n^{2/\beta} \ 
\varepsilon^{-2/\beta}.
$$
 After the minimization of the  right side over $ \varepsilon > 0 $ and $ (\alpha,\beta) $ 
we obtain 
$$
{\bf E}|\xi_n|^p \le 2 \cdot 5^p \  R^p(\gamma_n, p, K \times \psi(\cdot)), \ \
|\xi_n|_p \le 5 \sqrt{2} \ R(\gamma_n,p, K \times \psi(\cdot))
$$
and 
$$
||\xi_n||G(\nu) \le 5 \sqrt{2} \ \sup_{p \ge 2} V(\gamma_n,p, K \times \psi)/\nu(p) \to 0
$$
at $ n \to \infty $ by virtue of  condition (17). This completes the proof.\par
 Note that theorem 9 gives the concrete estimation $ ||S_n - S||G(\nu) $ in the term of
$ {\bf D}(S_n - S) = |S_n-S|_2^2: \ ||S_n - S||G(\nu) \le  $

$$
\le 5 \sqrt{2} \ \sup_{p \ge 2} \frac{V \left(|S_n - S|_2,p,[\sup_n ||S_n||G(\psi)] 
\times \psi(\cdot) \right) }{ \nu(p)}.
$$

 By virtue of Doob's inequality we can obtain analogous {\it maximal} inequality for the 
values $ \tau_n = || \ \sup_{l \ge n} |S_l - S| \ || G(\nu): $ \par 

$$
\tau_n \le 10 \sqrt{2} \sup_{p \ge 2} \frac{V(|S_n-S|_2,p,[\sup_n||S_n||G(\psi)]\times \psi(\cdot))}
{\nu(p)}.
$$

 For example  suppose that for some $ m > 0 \ \sup_n ||S_n||_m = C_1 < \infty. $ Let $ \Delta =
const > 0. $ Let us denote  $ \gamma_n = |S_n - S|_2. $ We obtain:

$$
|| \ S_n - S \ ||_{m/(m \Delta+1)} \le C(\Delta) \ |\log \gamma_n|^{-\Delta},
$$
and, moreover,

$$
|| \sup_{l \ge n}|S_l - S| \ ||_{m/(m \Delta + 1)} \le C(\Delta) \ |\log \gamma_n|^{-\Delta}.
$$

{\bf Corollary 1.} Since 
$$
V(\delta,p,\psi) \le C \ \delta^{1/(p+1)} \ \psi^{p/(p+1)}(2p),
$$
we obtain the following sufficient  condition for $ \nu(\cdot): $ if
$$
\lim_{\delta \to 0+} \sup_{p \ge 2} \delta^{1/(p+1)} \psi^{p/(p+1)}(2p)/\nu(p) = 0, \eqno(20)
$$
then for any martingale $ (S_n,F_n) $ such that 
$\sup_n ||S_n || G(\psi) < \infty $ it follows  $ ||S_n - S|| G(\nu) \to 0 $ at
$ n \to \infty. $\par
{\bf Corollary 2.} If $ \psi(\cdot) \in \Delta_2, $  then the condition 
$$
\lim_{p \to \infty} \psi(p) /\nu(p) = 0, \eqno(21)
$$ 
or, equally, $ N_{\nu}(\cdot) << N_{\psi}(\cdot) $  is sufficient for the  
implication (16). \par

 For instance, if $ \psi(p)= \psi_{m,r}(p) =  p^{1/m} \ \log^r p, \ m>0, r \in R^1, $ then 
$ \nu(p) $ may be, for example,
$$
\nu(p) =  p^{\Delta + 1/m} \ \log^r p, \ \  \nu(p) =  p^{1/m} \log^{\Delta+r} p, \ 
$$

$$
\nu(p) =  p^{1/m} \ \log^r p \ \log^{\Delta}(2+\log p), \ \ \Delta = const > 0.
$$
etc. In particular, if for some martingale $ (S_n,F_n) $ there 
exists $ m > 0 $ such that $ \sup_n ||S_n||_m < \infty, $ then $ \ \forall \Delta \in (0,m) $
$$
\lim_{n \to \infty} ||S_n - S||_{m-\Delta} = 0.
$$
 Further we obtain the necessary and sufficient conditions for implication (16) 
in the case if $ \psi \in EL_m. $ \\
 {\bf Theorem 10.} {\it Let  $ \psi \in EL_m $ for some $ m > 0; $ let $ \nu(\cdot) $
 be some function belonging $ \Psi. $ The following implication is true: \\
\vspace{3mm}
$ [\forall \ G(\psi) $ - bounded martingale $ (S_n,F_n):$

$$
 \sup_n ||S_n||G(\psi) < \infty \ \Rightarrow \eqno(22)
$$

$$
\ \Rightarrow \lim_n ||S_n - S||G(\nu) = 0] \eqno(23)
$$

if and only if}

$$
\lim_{p \to \infty} \psi(p)/\nu(p) = 0, \eqno(24)
$$
{\it or } $ N_{\nu} << N_{\psi}. $ \par
 {\bf Proof.} Sufficientness it follows immediately from theorem 7 and 
corollary 2. In order to prove necessarity of  condition (24) 
we must prove that if condition (24) is not satisfied then there exists a martingale 
$ (S_n,F_n) $ (which may be defined on {\it some} probability space)
such that $ \sup_n||S_n||L(N) < \infty $ but $ S_n $ does not converge
in the $ L(M) $ norm. Here $ N = N_{\psi}, M = N_{\nu}. $ Put for some  $ L_0 \in SL $ 
and $  B = const \in (0.5; 1)  $
$$
S_n = \sum_{k=2}^n k^{-B} L_0(k) \epsilon(k), \ \ S = \sum_{k=2}^{\infty} k^{-B} L_0(k) 
\epsilon(k), 
$$
where  $ \{\epsilon(k) \} $ is again  Rademacher sequence, $ F_n = \sigma \{\epsilon(i),
i \le n. $ We choose $ B $ and $ L_0(\cdot) $
 in the case $ m > 2 $ such that 
$$
B = m/(m+1), \ L_0^{-1/(1-B)} \left(u^{1/(1-B)} \right) = L(u),
$$
where $ N_{\psi}(u) = \exp \left(u^m L(u) \right), \ u \ge 2. $ \par
  From theorem 5 follows that 
$$
S \in G(\psi_N) \setminus G^0(\psi_N). 
$$
 Assume converse  for condition (24), i.e. that
$$
\overline{\lim}_{p \to \infty} \psi(p)/\nu(p) > 0.
$$
 Since $ ||S_n - S||G(\nu) \to 0 $ as $ n \to \infty $ and since $ S_n $ is bounded:
$$
vrai sup_{\omega \in \Omega} |S_n| \le \sum_{k=2}^n k^{-B} L_0(k),
$$
 the random variable $  S $ belongs to the space $ G^0(\nu).$ Therefore (see theorem 2)
$$
\lim_{p \to \infty} |S|_p/\nu(p) = 0.
$$
 But according to our conditions  it follows that 

$$
\overline{\lim}_{p \to \infty} |S|_p/\nu(p) > 0.
$$
 This contradictions  proves theorem 10, but only in the case $ m > 2. $  In order 
to prove our statement for the values $ m \in (1,2], $ we consider a {\it new} martingale 
$$
S^{(2)}_n = \sum \sum_{i,j = 1,2,...,n; i \ne j} i^{-B}L_0(i) \ j^{-B} L_0(j)\ \epsilon(i,1) 
\ \epsilon(j,2),
$$

$$
S^{(2)} = a.e. \  \lim_{n \to \infty} S^{(2)}_n 
$$
  with correspondence $ \sigma \ - $ flow $ F_n = \sigma \{ \epsilon(i,1), \epsilon(j,2); \ 
i,j \le n \}, $ where $ \epsilon(i,s) $ are independent sequences (over $  s = 1,2,...$ ) of 
 Rademacher series, $ B $ again belongs to the interval $ (0.5; 1). $ It follows from theorem 7
and the representation $ S_n^{(2)} = $ 
$$
= \sum_{i=1}^n i^{-B}L_0(i) \epsilon(i,1) \times \sum_{j=1}^n 
j^{-B} L_0(j) \epsilon(j,2)- \sum_{i=1}^n i^{-2B}L_0^2(i) \epsilon(i,1) \epsilon(i,2)
$$
that for all $ x \ge 2 $

$$
\exp \left(- C_1x^{m/2} L(\sqrt{x}) \right) \le  {\bf P} (|S^{(2)}| > x) \le \exp \left(-C_1 
x^{m/2} L(\sqrt{x}) \right).
$$

 As above we conclude that if 
$$ 
\overline{\lim}_{p \to \infty} \psi_N(p)/\nu_M(p) > 0,
$$
then  $ S^{(2)} \in EOS(M_{m/2,L(\sqrt{\cdot})} \setminus EOS^0(M_{m/2,L(\sqrt{\cdot})}). $ 
Therefore 
$$
\lim_{n \to \infty} ||S^{(2)}_n - S^{(2)}||L(N_{m/2,L(\sqrt{\cdot})}) \ne 0
$$
 and now $ m/2 \in (1;\infty). $  Analogously the case $ m \in (0.5; \infty) $ etc. may be 
considered. \par

 {\it Remark 4.} The same result as in theorem 10 is also true in the case $ V(Z,\beta) $
spaces. Namely, for  all martingales $ (S_n,F_n) $ from condition

$$
\sup_n ||S_n||G \left(\psi^{(Z,\beta)} \right) < \infty
$$
 follows that for some $ \nu(\cdot) \in \Psi $

$$
\lim_{n \to \infty} ||S_n - S||G(\nu) = 0
$$
if and only if $ \lim_{p \to \infty} \psi^{(Z,\beta)}(p)/\nu(p) = 0. $ \par
 
 The conclusion "if"  follows from theorem 9, the counterexample in the spirit of 
theorem 10 may be constructed by formula 
$$
Y_n = \sum_{d=1}^{\infty} C(d) \ Y(d,n), \ \ Y = a.e. \lim_{n \to \infty} Y_n.
$$
where 

$$
Y(d,n) = \sum \sum \ldots \sum_{k_i = 1,2,...,n; k_i \ne k_j, i \ne j} \frac{\prod_{i=1}^d 
\epsilon(k_i,i,d)}{ \prod_{i=1}^d k_i^B },
$$

$ \epsilon(k,s,d) $ are  a family (over $ (s,d)) $ 
 of independent Rademacher sequences,  $ B \in (1/2;1), 
m = 1/(1-B) \in (2, \infty),  \ C(d) $ is some sequence of constant (we will choose $ C(d) $ 
further) and correspondence $ \sigma $ - flow $ \{ F_n \} $ is a natural; hence the pair 
$ (Y_n, F_n) $ is  a martingale. \par 
 From the proof of theorem 10 we obtain that $ \forall u \ge 2 $

$$
\exp \left(-C_2^d \ u^{m/d} \right) \le {\bf P}(|Y(d,n)| > u) \le \exp \left(-C_1^d \ u^{m/d}
\right).
$$
 By virtue of theorem 1 we  receive the bide - side moment estimations:

$$
 C_3 \ p^{d/m} \le  |Y(n,d)|_p \le C_4 \ p^{d/m}, \ p \ge 2. 
$$

 Now let us choose $ C(d): \ C(d) = d^{-d \gamma}, \ \gamma = const > 0. $ We obtain from 
the triangular inequality:

$$
|Y_n|_p \le \sum_{d=1}^{\infty} C^d \left(p^{1/m} \right)^d \  d^{-d\gamma} \le 
 \exp \left(C_5(\gamma,m) p^{1/(m \ \gamma)}  \right) \  \eqno(25)
$$
(upper bound). Now we obtain the low bound for $ |Y_n|_p. $ Since the martingales $ \{Y(d,n) \} $
are independent  $ |Y_n|_p \ge \sup_d \ |Y(d,n)|_p \ge $

$$
 \ge \sup_d C^d \left( p^{1/m}  \right)^d \ d^{-d\gamma} \ge
 \exp \left( C_4(\gamma,m) p^{1/(m \ \gamma)} \right). \eqno(26)
$$

 From estimations (25), (26) follows that 

$$
Y \in V(Z,\beta) \setminus V^0(Z,\beta), \ \beta = 1/(m\gamma), \ \exists Z \in (0,\infty).
$$

 Therefore, if 
$$
\overline{\lim}_{p \to \infty} \psi^{(Z,\beta)}(p)/\nu(p) > 0
$$ 

then
$$
\lim_{n \to \infty} ||Y_n - Y||G({\nu}) \ne 0. 
$$

\vspace{3mm}

{\bf Aknowledgements.} I am very grateful to prof. V. Fonf,  M.Lin, B.Rubinstein for 
useful  support of these investigations. I am also very grateful to 
prof. N. Kasinova for her  consultations. \par
 This work  was partially supported by the Israel Ministry of Absorbtion of Israel.\par

\vspace{3mm}

\newpage
\begin{center}
{\sc REFERENCES} \\
\end{center}
\vspace{2mm}
1. H.Aimar, E.Harbour and B. Iaffel.   Boundedness of Convolution Operators with 
smooth Kernels on Orlicz Spaces. {\it Studia Math.,} 151(3),  (2002), 195 - 206.\\ 
2. D.R. Bagdasarov, E.I. Ostrovsky.  Reversion of Chebyshev's Inequality.
{\it Probab. Theory Appl.,} v.40  $ N^o $ 4, 737 - 742.\\
3. V.V.Buldygin, D.I.Mushtary, E.I.Ostrovsky, M.I.Puchalsky.  New Trends in 
Probability Theory and Statistic. MOKSLAS, 1992, Amsterdam, New York.\\
4. V.V.Buldygin, Ju.V.Kozachenko.  Metric Characterization of Random Variables 
and Random Processes. AMS, 678, 2000, Providence, RI.\\
5. R.E.Edwards.  Fourier Series. A Modern Introduction . Springer Verlag, 1982, 
 Berlin, Heidelberg, Hong Kong.\\
6. Y.Giga, H.Sohr.  Abstract $ L^p $ Estimations for Cauchy Problem with 
Applications to the Navier - Stokes Equations in Exterior Domains. {\it J. Funct. 
Anal., } 102, $ N^o $ 1, (1991), 72 - 94.\\
7. P.Hall, C.C.Heyde.  Martingale Limit Theory and its Applications. Academic Press,
1979,  New York, London, Toronto, Sydney, San Francisco.\\
8. M.A. Krasnoselsky, Ya.B.Routisky.  Convex Functions and Orlicz Spaces.
P.Noordhoff Ltd, 1961, Groningen.\\
9.A.Kufner, O.John, S.Fuchik.  Function Spaces. Academia, Prague
 and Noordhoff Ltd,  1979, Groningen.\\
10. J.Lindenstraus, L.Tsafriri.  Classical Banach Spaces. Springer Verlag, 
1977, Berlin - Heidelberg, New York.\\
11. J.Neveu. Discrete - Parameter Martingales. North - Holland Publishing Company,
1975, Amsterdam - Oxford - New York.\\
12. E.I.Ostrovsky. Exponential Estimations for the Random Fields. OINPE, 
1999,  Obninsk (in Russian).\\
13. G.Peshkir. Maximal Inequalities of Kahane - Khintchin type in Orlicz Spaces.
Preprint Series $ N^o $ 33, Institute of Mathematics, 1992, University of Aarhus (Danemark).\\
14 S.K.Pichorides. On the best values of the constant in the theorem of M.Riesz, 
Zygmund and Kolmogorov. {\it Studia Math.,} 44, (1972), 165 - 179.\\
15. M.M.Rao, Z.D.Ren.  Theory of Orlicz Spaces. Marcel Dekker Inc., 1991, 
New York, Basel, Hong Kong.\\
16. M.M.Rao, Z.D. Ren. Applications of Orlicz Spaces. Marcel Dekker Inc., 2002,
 New York, Basel.\\
17. E.Seneta.  Regularly Varying Functions. Springer Verlag, 1976, Berlin, 
Heidelberg, New York.\\
18. E.M.Stein.  Singular Integrals and Differentiability Properties of Functions.
Princeton Univ. Press, 1970, Princeton, New York.\\

\vspace{3mm}
 ISRAEL, Beer - Sheva, Ben Gurion University,\\
84105, Ben Gurion street, 4; Postal Box 61.\\
Dept. of Math. and Comp. Science.\\
e - mail: Galaostr@cs.bgu.ac.il.\\

\newpage

\begin{center}
{\bf Exponential Orlicz Spaces: New Norms and  Applications.}\\ 
\vspace{3mm}
{\bf E.I.Ostrovsky.} \\
\vspace{2mm}
Abstract.
\end{center}

 Let $ (\Omega, F, {\bf P}) $ be  probability space,  $ N = N(u) $ be an exponential $ N \ - $
Orlicz function. We introduce in the Orlicz space $ (\Omega, F, {\bf P}, N) $ 
other norms  which are equivalent to the classical Orlicz norm 
$ ||\cdot||L(N),  $ for example, by means of all moments:
$$
||\eta||G(\psi) = \sup_{p \ge 2} |\eta|_p/\psi(p), \ |\eta|_p = {\bf E}^{1/p} |\eta|^p,
$$
 and show  convenience of their  applications in the  theory of Orlicz spaces,
in the operator theory, in  the  theory of  Fourier series  and in the theory of martingales. \par
 For instance, let $ (S_n,F_n) $ be a martingale over some probability space  
and $ N = N(u) $ be some exponential Orlicz function  such that 
$$
\sup_n ||S_n||L(N) < \infty, \ N \notin \Delta_2.
$$
 where $ ||\cdot|| $ is  the classical Luxemburg norm. We study all new $ N  $ functions
$ M = M(u) $ such that for all martingales $ (S_n, F_n) $
$$
\lim_{n \to \infty} || S_n  - S ||L(M)  = 0, \ S = a.e. \ \lim_{n \to \infty} S_n.
$$
and show that in general case of exponential Orlicz function $ N(u) $ the condition $ M(\cdot) 
 << N(\cdot) $ is necessary and sufficient for this implication.\par

 {\it References:} 18 publications. 

\end{document}